\documentclass{article}

\title{Transformations of algebraic Gauss
hypergeometric functions}

\author{Raimundas Vid\=unas
\\ \em Kobe University}

\newtheorem{theorem}{Theorem}[section]

\newtheorem{algorithm}[theorem]{Algorithm}
\newtheorem{example}[theorem]{Example}
\newtheorem{definition}[theorem]{Definition}

\usepackage{latexsym}

\newcommand{\hpg}[5]{{}_{#1}\mbox{\rm F}_{\!#2}\!
  \left(\left.{#3 \atop #4}\right| #5 \right) }

\newcommand{\hpgo}[2]{{}_{#1}\mbox{\rm F}_{\!#2}}

\newcommand{\qed}{\hfill $\Box$\\}
\newcommand{\algstep}{\vspace{5pt}\\}
\newcommand{\equal}{&\!\!\!=\!\!\!&}
\newcommand{\CC}{\mbox{\bf C}}

\newcommand{\PP}{\mbox{\bf P}}

\newcommand{\DD}{\mbox{\bf D}}

\begin{document}

\maketitle

\begin{abstract} A celebrated theorem of Klein implies that
any hypergeometric differential equation with algebraic solutions
is a pull-back of one of the few standard hypergeometric equations
with algebraic solutions. The most interesting cases are
hypergeometric equations with tetrahedral, octahedral or
icosahedral monodromy groups. We give an algorithm for computing
Klein's pull-back coverings in these cases, based on certain explicit
expressions (Darboux evaluations) of algebraic hypergeometric functions.
The explicit expressions can be computed from a data base (covering the Schwarz table)
and using  contiguous relations.
Klein's pull-back transformations also induce algebraic transformations between
hypergeometric solutions and a standard hypergeometric function 
with the same finite monodromy group.
\end{abstract}

\section{Introduction} \label{alghypers}

Algebraic Gauss hypergeometric functions is a very classical
object of research in mathematics. They were studied by many
authors \cite{schwarz72,fuchs,brioschi,
klein78,pepin,boulang,katz72,baltadw,singulm2,putulm}. Recall that
Gauss hypergeometric function $\hpg{2}{1}{\!A,\;B\,}{C}{Z\,}$
satisfies the {\em hypergeometric differential equation}:
\begin{equation} \label{hpgde}
Z\,(1-Z)\,\frac{d^2Y(Z)}{dZ^2}+
\big(C-(A\!+\!B\!+\!1)\,Z\big)\,\frac{dY(Z)}{dZ}-A\,B\;Y(Z)=0.
\end{equation}
This is a Fuchsian equation 
with three regular singular points $Z=0,1$ and $\infty$ on
$\PP^1$. The local exponent differences at these points are (up to
a sign) $1-C$, $C-A-B$ and $A-B$ respectively. We denote equation
(\ref{hpgde}) by $H(1-C,C-A-B,A-B)$. With this notation, equation
the $H(e_0,e_1,e_\infty)$ can be written in the following form:
\begin{equation} \label{hpgdee}
\frac{d^2Y(Z)}{dZ^2}+ \left( \frac{1-e_0}{Z}+\frac{1-e_1}{Z-1}
\right)\,\frac{dY(Z)}{dZ}+\frac{(1-e_0-e_1)^2-e_\infty^2}
{4\,Z\,(Z-1)}\;Y(Z)=0.
\end{equation}
If we permute the parameters $e_0,e_1,e_\infty$ or multiply some
of them by $-1$, we get hypergeometric equations related by
well-known fractional-linear transformations, see \cite{specfaar}.
In general, there are 24 hypergeometric equations related in this
way, and they share the same (up to fractional-linear transformations of
the argument and radical factors) 24 hypergeometric Kummer's solutions.

As is known (see \cite{beukers} for example), a Fuchsian equation
only has algebraic solutions if and only if its {\em monodromy
group} is
finite. 
The following hypergeometric equations (and their
fractional-linear transformations just mentioned) have this
property and are called {\em standard hypergeometric equations}
with algebraic solutions:
\begin{itemize}
\item $H(1,1/n,1/n)$, where $n$ is a positive integer. The
hypergeometric equation degenerates to a Fuchsian equation with
two singular points. Its monodromy group is the cyclic group with
$n$ elements.
\item $H(1/2,1/2,1/n)$, where $n$ is an integer, $n\ge 2$. The
monodromy group of this equation is the dihedral group with $2n$
elements.
\item $H(1/2,1/3,1/3)$. The monodromy 
group is the tetrahedral group, isomorphic to $A_4$.
\item $H(1/2,1/3,1/4)$. The monodromy 
group is the octahedral group, isomorphic to $S_4$.
\item $H(1/2,1/3,1/5)$. The monodromy 
group is the icosahedral group, isomorphic to $A_5$.
\end{itemize}
The celebrated theorem of Klein \cite{klein77} states that if a
second order linear homogeneous equation only has algebraic
solutions, then it 
is a {\em pull-back transformation}
of a standard hypergeometric equation from the list above.
If the Fuchsian equation has coefficients in $\CC(X)$, such a
pull-back transformation changes the variable $Z$ in
(\ref{hpgdee}) to a rational function $\varphi(X)$. In geometric
terms, we have a finite covering $\varphi:\PP^1_X\to\PP^1_Z$
between two projective lines, and we {\em pull-back} the standard
hypergeometric equation from $\PP^1_Z$ onto $\PP^1_X$.

Particularly, the theorem of Klein implies that if a
hypergeometric equation $H_1$ only has algebraic solutions, then
it is a pull-back transformation of a standard hypergeometric
equation $H_0$ with algebraic solutions. The proof of Klein 
implies that the monodromy group of $H_1$ is either
cyclic, or dihedral, or $A_4$, $S_4$ or $A_5$, and that $H_0$ can
be chosen to be a standard hypergeometric equation with the same
monodromy group. The main purpose of this paper is to give an
algorithm for computing the pull-back transformations implied by
the theorem of Klein. We consider here the most interesting cases
of hypergeometric equations with tetrahedral, octahedral or
icosahedral monodromy groups. Hypergeometric equations with
dihedral 
monodromy groups are considered in \cite{tdihedral}.

Hypergeometric equations with algebraic solutions 
were first classified by Schwarz in \cite{schwarz72}. Disregarding
hypergeometric equations with a cyclic monodromy group, Schwarz
gave a list of fifteen types of these hypergeometric equations.
One type consists of hypergeometric equations with a dihedral
monodromy group. The other types are represented by the following
hypergeometric equations:
\begin{itemize}
\item $H(1/2,1/3,1/3)$, $H(1/3,1/3,2/3)$. The monodromy
group is the tetrahedral group.
\item $H(1/2,1/3,1/4)$, $H(2/3,1/4,1/4)$. The monodromy
group is the octahedral group.
\item $H(1/2,1/3,1/5)$, $H(1/2,1/3,2/5)$, $H(1/2,1/5,2/5)$,
$H(1/3,1/3,2/5)$, $H(1/3,2/3,1/5)$, $H(2/3,1/5,1/5)$,
$H(1/3,2/5,3/5)$, $H(1/3,1/5,3/5)$, $H(1/5,1/5,4/5)$,
$H(2/5,2/5,2/5)$. The monodromy 
group is the icosahedral group.
\end{itemize}
We refer to Schwarz type of hypergeometric equations with
algebraic solutions by the triple of the parameters
$e_0,e_1,e_\infty$ 
of these representative equations. (Usually, Schwarz type is a
roman numeral from I to XV.) Hypergeometric equations of the same
Schwarz type are characterized by the property that their
hypergeometric solutions are {\em contiguous} to hypergeometric
solutions of the representative hypergeometric equation of that
type. 
For example, the parameters $e_0,e_1,e_\infty$ of hypergeometric
equations of the Schwarz type $(1/3,1/3,2/3)$ can be
characterized as follows: they are rational numbers, their
denominators are equal to 3, and the sum of their numerators is
even. See \cite[Subsection 5.6]{dalggaus} for more details.

Recently,  M. van Hoiej, J.-A.~Weil and M.~Berkenbosch 
\cite{WeHuBe}, \cite{kleinvhw} 
developed an algorithm for computing Klein's pull-back transformations
for second order linear differential equations, based on computation of
monodromy (or differential Galois group) semi-invariants,
via explicit solution of suitable symmetric tensor powers of those differential equations.
extension of Klein's theorem to third order
linear differential equations equations with a finite monodromy group is possible \cite{maintphd}.

The algorithm in this paper is devised for hypergeometric equations with algebraic solutions.
Instead of finding semi-invariants by solving symmetric tensor power equations,
we propose to use contiguous relations and a data base of explicit expressions for simplest
hypergeometric solutions. This approach should be more effective
when local exponent differences of hypergeometric equations are large,
since solving symmetric power equations apparently involves
a linear algebra problem whose size is proportional to the local exponent differences,
while deriving a contiguous relation can be be of logarithmic complexity \cite{contiguous}.
On the other hand, direct transformation of explicit expressions by Gr\"obner basis methods
appears border-line cumbersome for differential equations with the icosahedral monodromy group.
But these transformations can be apparently streamlined in each of the icosahedral Schwartz cases using a normal series of the solvable icosahedral group, as elimination degrees are universally bounded. This approach simplifies greatly
for algebraic dihedral hypergeometric functions, as demonstrated in \cite{tdihedral}.

The input of our algorithm is a hypergeometric equation
with a  tetrahedral, octahedral or icosahedral monodromy group. 
The output is a finite covering $\varphi:\PP^1\to\PP^1$ which pull-backs the corresponding
standard hypergeometric equation to the given one. With a little
more work one may start with an algebraic hypergeometric Gauss
function, and produce an algebraic transformation which relates
the given hypergeometric function to hypergeometric solutions of
the corresponding standard hypergeometric equation.
The algorithm is based on the following three ingredients:
\begin{itemize}
\item {\em Pull-back transformation of a hypergeometric equation
to a Fuchsian equation $($not necessarily on $\PP^1)$ with a
cyclic monodromy group.} Geometrically, we pull-back the
hypergeometric equation onto its {\em Darboux curve} by a
{\em Darboux covering}; see 
Section \ref{someprelim}. 
The algebraic hypergeometric solutions are transformed to {\em
radical functions} (i.e., products of powers of rational functions
on the Darboux curve). This transformation allows us to write
algebraic Gauss functions in a convenient explicit form.
General theory and computations of Darboux coverings and evaluations
is presented in \cite{dalggaus}.
\item {\em Contiguous relations of Gauss hypergeometric functions.}
With a pre-computed data basis of Darboux coverings 
and explicit Darboux evaluations for simplest algebraic hypergeometric functions of
each Schwarz type at hand, we use contiguous relations to compute any
algebraic algebraic function in that explicit form. See
\cite[Subsection 5.6]{dalggaus}. 
\item {\em Schwarz map}. 
As a function, this is a quotient of two solutions of the same
hypergeometric equation. If $s_0$ is a Schwarz map for a standard
hypergeometric equations with algebraic solutions, then $s_0^{-1}$
is a finite Galois covering $\PP^1\to\PP^1$ on which the monodromy
group acts as the Galois group. Schwarz map plays an essential
role in the proof of Klein \cite{klein77}; see Section
\ref{someprelim}. Klein's pull-back covering is identical to a Schwarz map of
the given hypergeometric equation composed with the rational
function $s_0^{-1}$ of the corresponding standard equation. We
make this construction very explicit using the previous two
ingredients.
\end{itemize}
Our algorithm is presented in Section \ref{kleinmorphs}. Sections \ref{sec:examples}
and \ref{sec:ixamples} demonstrate the algorithm by several examples. In Section
\ref{database} we present the mentioned data base of Darboux
coverings and explicit hypergeometric evaluations.

\section{Some preliminaries} \label{someprelim}

Here we bring forward the most important notions for
our algorithm. A more comprehensive presentation of related theory
(with further references), and computation of the data base in Section \ref{database}
are given in the supplementing paper \cite{dalggaus}.

First we briefly sketch the Klein's proof of his theorem in
\cite{klein77} and \cite[Section V.3]{icosaklein}. The theorem
states that if the differential equation $Y''+p\,Y'+q=0$ with
$p,q\in\CC(X)$ only has algebraic solutions, it is a pull-back
transformation of one of the standard hypergeometric equations.
Let $s\in\CC(X)$ denote a quotient of two distinct solutions of
this equations. Then $s$ satisfies the equation
\begin{equation} \label{schwder}
\frac{s'''}{s'} -
\frac32\left(\frac{s''}{s'}\right)^2=2q-p'-\frac12\,p^2.
\end{equation}
Incidentally, the expression on the left-hand side is the {\em
Schwarzian derivative} $\{s,X\}$. Solutions of this equation are
precisely the fractional-linear transformations $(as+b)/(cs+d)$ of
$s$. Therefore, analytic continuation of $s$ along a closed loop
(away from singularities of the given equation) changes $s$ by a
fractional-linear transformation. All closed paths with a common
fixed point give a group of fractional-linear transformations that
acts on the branches of $s$. This is the {\em projective monodromy
group} of the equation; see \cite[Section 5.4]{dalggaus} for
example. Since $s$ is an algebraic function, this group must be
finite. Finite groups of fractional-linear transformations are
well-known: such a group is either cyclic, or dihedral, or $A_4$,
$S_4$ or $A_5$. Functions in $\CC(s)$ that are invariant under the
transformation group are rational functions of $X$. When a
suitable invariant function is equal to $X$, the original equation
$Y''+pY'+q=0$ is a standard hypergeometric equation with the same
projective monodromy group. In general, the same invariant
function (as a rational function of $X$) gives the required
pull-back transformation from the same standard hypergeometric
equation.

If the equation $Y''+pY'+q=0$ is a hypergeometric equation, we
interpret the construction of Klein as follows. The quotient $s$
of two solutions is a {\em Schwarz map} for the hypergeometric
equation, see \cite[Section 5.5]{dalggaus} for example. If $s_0$
is the quotient of two corresponding solutions of the standard
hypergeometric equation with the same monodromy group, then
$s_0^{-1}$ is the ``suitable" invariant function that cancels the monodromy of $s$.
Klein's pull-back covering is the composition $s_0^{-1}\circ s$.

A way to obtain Klein's pull-back covering is to compute the quotient $s$
of two convenient hypergeometric solutions 
as a root of an algebraic equation, and then use elimination
techniques to compute the composition $s_0^{-1}\circ s$ as a
rational function of $X$. 
We particularly propose to use a convenient expression for the function 
$s$ by invoking {\em Darboux curves} and {\em Darboux coverings}, 
so to make computation of $\psi(s)$ algebraically straightforward. 
Specifically, we propose to pull-back the given differential equation to its Darboux curve,
and express its hypergeometric solutions and the quotient $s$ of
two solutions as radical functions on the Darboux curve. This
is a very convenient and explicit expression for the function $s$. 

In the remainder of this Section we present and comment our
working definition of Darboux curves and Darboux coverings 
and sketch their application to our algorithm.
We refer to \cite[Section 3]{dalggaus} for a thorough discussion of Darboux curves and Darboux coverings, and to \cite[Section 5.3]{dalggaus} for our specifications of
pull-back transformations.
\begin{definition} \rm
Let $E$ denote a linear homogeneous differential equation with
coefficients in $\CC(X)$. Let $\DD$ be an algebraic curve, and let
$\CC(\DD)$ be the field of rational functions on $\DD$. Suppose
that $\varphi:\DD\to\PP^1$ is a finite covering. We say that $\DD$
is a {\em Darboux curve} for $E$, and that $\varphi$ is a {\em
Darboux covering} for $E$, if a pull-back transformation of $E$
onto $D$ with respect to $\varphi$ has a solution $Y$ which
satisfies $Y'=u\,Y$ for some $u\in\CC(\DD)$, and the algebraic
degree of $u$ over $\CC(X)$ is equal to the degree of $\varphi$.
\end{definition}

In general, a pull-back transformation of a differential equation allows projective transformation
$y(x)\mapsto\theta(x)y(x)$ of solutions by a radical function $\theta(x)$. 
A general pull-back transformation has the form
\begin{equation}
z\mapsto \varphi(x), \qquad  y(z)\mapsto \theta(x)\,y(\varphi(x)).  
\end{equation}

When speaking about Darboux curves, it is more convenient to refer
to the {\em differential Galois group} of a differential equation
rather than to the monodromy group. In the case of hypergeometric
equations with algebraic solutions only, these two groups
coincide. 
Particularly, the Piccard-Vessiot extension (where the solutions
``live") is an algebraic Galois extension, it is generated over
$\CC(X)$ by a solution of the hypergeometric equation, and its
usual Galois group is the same monodromy group. For basic notions
of differential Galois theory and further references (in particular, to proofs of
mentioned facts) 
we refer to \cite[Section 5.2]{dalggaus}.

Consider a hypergeometric equation $H_1$ with algebraic solutions
only. Let $K\supset \CC(X)$ be its Piccard-Vessiot extension, and
let $G$ denote its differential Galois group. If $\DD$ is a
Darboux curve for $H_1$, then its function field $\CC(\DD)$ is an
intermediate field: $K\supset\CC(\DD)\supset\CC(X)$. Besides,
$K\cong\CC(X,Y)$ and $\CC(\DD)\cong\CC(X,u)$ for $u=Y'/Y$ and for
a suitable solution $Y$ of $H_1$. The extension $K\supset\CC(\DD)$
must be a finite subgroup of the multiplicative group
$\mbox{GL}(1,\CC)$, so it must be a finite cyclic group.
Conversely, each cyclic subgroup of $G$ defines a Darboux curve.
In our algorithm we use Darboux curves which correspond to cyclic
subgroups of maximal order. These Darboux curves have the smallest
genus, and the corresponding Darboux coverings have minimal
degree; see \cite[Table 1]{dalggaus}.

From now on, let $X$ denote the indeterminate of a hypergeometric
equation with algebraic solutions, and let $Z$ denote the
indeterminate of a standard hypergeometric equation with the same
differential Galois group $G$. By \cite[Lemma 3.3\em (i)]{dalggaus}, 
all Darboux curves for the standard hypergeometric equation are rational; 
we denote a rational parameter of such a Darboux curve by $z$. 
In the current context, let $H_0$ denote the standard hypergeometric equation for $H_1$,
and consider a Darboux covering $\varphi_0:\PP^1_z\to\PP^1_Z$ for
$H_0$ of minimal degree. Let $\psi:\PP^1_X\to\PP^1_Z$ be the Klein
covering which pull-backs $H_1$ to $H_0$, and let
$\varphi_1:\DD\to\PP^1_X$ be the Darboux covering for $H_1$ of
minimal degree. By \cite[Lemma 3.5]{dalggaus}, the Darboux curve
$\DD$ for $H_1$ is the fiber product of
$\varphi_0:\PP^1_z\to\PP^1_Z$ and
$\psi:\PP^1_X\to\PP^1_Z$, and $\varphi_1$ 
is the projection onto $\PP^1_X$. We have the following diagram:
\begin{equation} \label{di:darboux}
\begin{picture}(270,95)(0,3)
\put(144,-1){$\PP^1_Z$} \put(114,40){\vector(1,-1){30}}
\put(186,40){\vector(-1,-1){30}} \put(106,45){$\PP^1_z$}
\put(187,45){$\PP^1_X$} \put(145,86){\vector(-1,-1){30}}
\put(157,86){\vector(1,-1){30}} \put(147,90){$\DD$}
\put(71,90){$\PP^1_s$} \put(75,86){\vector(1,-1){30}}
\put(119,21){$\varphi_0$}\put(174,21){$\psi$}
\put(80,62){$\gamma$} 
\put(172,75){$\varphi_1$}
\end{picture}
\end{equation}
Here the composition $\phi_0\circ\gamma:\PP^1_s\to\PP^1_Z$ is the
inverse Schwarz map $s_0^{-1}$ for $H_0$, and formally, a Darboux
covering for $H_0$ of maximal degree. Besides, $\CC(s)\supset\CC(Z)$
is the Piccard-Vessiot extension for the standard hypergeometric
equation. By \cite[Lemma 3.3]{dalggaus}, the map
$\gamma:\PP^1_s\to\PP^1_z$ is a Galois covering, and its Galois
group is a maximal cyclic subgroup of $G$.

The role of Darboux coverings in our main Algorithm \ref{kleinalg}
can be viewed in two ways. On one hand, we can identify the
hypergeometric solutions of $H_1$ and $H_0$ up to unknown
scalar  factor; 
see our examples in Section \ref{sec:examples} below. 
The scalar factor can be determined by evaluating the hypergeometric 
functions on some corresponding points on the Darboux curves $\DD\to\PP_z^1$.
Once we choose a point $P$ on $\PP_z^1$, there are (a priori) finitely many choices on $\DD$
in the same fiber of $\varphi_1\circ\psi$ above $\varphi_0(P)$ for this identification. 
A further detailed analysis of all Schwartz cases may define unique identifications here.
To identify a right pair of corresponding points, one can check the possibilities using 
a pair of compatible  automorphisms of the two Darboux coverings.

Let $f_1,f_2$ be two hypergeometric solutions of $H_1$ identified in this way with two
hypergeometric solutions $g_1,g_2$ (respectively) of $H_0$. The
pull-backs of $f_1,g_1$ (and of $f_2,g_2$) onto the respective
Darboux curves $\PP^1_z$, $\DD$ can be identified, once the scalar factor is determined.
It is important in  computations that those pull-back functions are radical functions
on the Darboux curves, so they have convenient explicit
expressions. The right identification of the quotients $f_1/f_2$ and $g_1/g_2$ determines
Klein's pull-back covering $\psi$. This identification of 
$f_1/f_2$ and $g_1/g_2$ on the Darboux curve $\DD$, and
explicit expressions for $\varphi_1$, $\varphi_0$ give enough
equations to determine an algebraic relation between $X$ and $Z$.
Klein's theorem ensures that factorization of this relation gives
the expression $Z=\psi(X)$ for Klein's pull-back covering.

Other interpretation is conceptually more straightforward. Note that $f_1/f_2$
and $s_0=g_1/g_2$ are Schwarz maps for $H_1$ and $H_0$
respectively. The construction of Klein implies that
$s_0^{-1}(f_1/f_2)$ is precisely Klein's covering, once the two quotients
are rightly identified. We use Darboux curves here to have a convenient explicit expression 
for the quotient $f_1/f_2$ as a function on $\DD$. We can make
standardized identifications and use the most handy Galois invariant functions
$s_0^{-1}$ for the tetrahedral, octahedral and icosahedral cases.
Klein's theorem ensures that in our expression of
$s_0^{-1}(f_1/f_2)$ we can eliminate the variables of $\CC(\DD)$
and get the required function $\psi(X)$. This interpretation gives
This is a convenient conceptual interpretation for our Algorithm \ref{kleinalg}.

The relation of Schwarz maps and Klein's pull-back covering is surely widely known. Our
use of Darboux coverings for explicit computations is perhaps new. It is
interesting to note that in \cite{ochiay} Ochiai and Yoshida compare the
Schwarz maps $s$ and $s_0$ by considering the dual composition $s\circ
s_0^{-1}$. They consider only hypergeometric equations of standard Schwarz
types, so their hypergeometric functions are contiguous to hypergeometric
solutions of the standard equations $H(1/2,1/3,1/3)$, $H(1/2,1/3,1/4)$,
$H(1/2,1/3,1/5)$. Then the composition $s\circ s_0^{-1}$ is a finite
covering $\PP^1\to\PP^1$ as well. The finite coverings $s_0^{-1}\circ s$,
$s\circ s_0^{-1}$ are related by fiber product construction via two copies
of $s_0^{-1}$:
\[
\begin{picture}(220,67)(8,38)
 \put(108,45){$\PP^1_x$}
\put(187,45){$\PP^1_z$} \put(145,86){\vector(-1,-1){30}}
\put(158,85){\vector(1,-1){29}} \put(147,90){$\PP^1_s$}
\put(71,90){$\PP^1_{s_0}$} \put(79,85){\vector(1,-1){29}}
\put(122,48){\vector(1,0){60}} \put(85,93){\vector(1,0){58}}
\put(95,74){$s_0^{-1}$} 
\put(100,97){$s\circ s_0^{-1}$} 
 \put(133,66){$s^{-1}$}
\put(134,38){$s_0^{-1}\circ s$} \put(174,73){$s_0^{-1}$}
\end{picture}
\]
Here on $\PP^1_x$ we consider a hypergeometric equation with
algebraic solutions, on $\PP^1_z$ --- the corresponding standard
equation, and on $\PP^1_{s_0}$ --- both hypergeometric equations. The
covering $s\circ s_0^{-1}$, as a rational function of $s_0$, can
be conveniently computed in a similar way, using the same
explicit expression of $s$ as a function on the Darboux curve.

\section{Computation of Klein's pull-back coverings}
\label{kleinmorphs}

Here we present an algorithm for computing Klein's coverings.  
As mentioned in the introduction, it is based on explicit evaluations of algebraic Gauss
hypergeometric functions in Section \ref{database} (which use
Darboux coverings); contiguous relations (by which we get similar
explicit expressions for any hypergeometric functions); and
identification of Schwarz maps  (to get an extra algebraic
relation between parameters on different curves).
\begin{algorithm} \rm \label{kleinalg} 
{\em Input:} A hypergeometric equation $H_1$ with
tetrahedral, octahedral or icosahedral 
monodromy group.
\algstep {\em Output:} A finite covering $\psi:\PP^1\to\PP^1$, which
pull-backs a standard hypergeometric equation $H_0$ with the same
monodromy group to $H_1$. 
\algstep {\bf Step 0}. As in the previous section, we denote the variables
of the differential equations $H_0$ and $H_1$ by $Z$ and $X$
respectively. The local exponent differences of $H_0$ are
$1/2,1/3,1/m$, where $m=3,4$ or $5$ for the
tetrahedral, octahedral or icosahedral monodromy groups, respectively. 
We assign the local exponent differences $1/2,1/3,1/m$ to the points $Z=1$,
$Z=\infty$ and $Z=0$, respectively.
\algstep {\bf Step 1.} We assume that the point $X=0$ of the given equation 
lies above the point $Z=0$ under Klein's pull-back covering.
The points $X=1$ and $X=\infty$ must be assigned so that
they would lie above the $Z$-points with the same denominators 
of local exponent differences. Particularly, if the given equation 
has a half-integer local exponent, we assume that the point $X=1$ has that local exponent,
and that it lies above $Z=1$. 
\algstep {\bf Step 2.} We use the Darboux coverings and evaluations
explicitly presented in Section \ref{database}. Let $z$ denote
a rational parameter for the Darboux curve for the standard equation $H_0$.
We denote coordinate functions of the Darboux curve for $H_1$ by $x$ and $\xi$. 
(If the Darboux curve is rational, then $\xi=0$ and we use the variable $x$ only.) 
Let $\phi_0(z)$, $\phi_1(x,\xi)$ denote the Darboux coverings for $H_0$, $H_1$ respectively. 
We have the algebraic relations $Z=\phi_0(z)$, $X=\phi_1(x,\xi)$ and
an algebraic relation between $x$, $\xi$ (which possibly $\xi=0$).
\algstep {\bf Step 3.} Consider the following solutions of $H_0$:
\begin{equation} \label{twoh1sols}
F_1=\hpg{2}{1}{\frac{1}{12}-\frac{1}{2m},\,\frac{5}{12}-\frac{1}{2m}}
{1-\frac{1}{m}}{Z}, \quad
F_2=Z^{1/m}\,\hpg{2}{1}{\frac{1}{12}+\frac{1}{2m},\,\frac{5}{12}+\frac{1}{2m}}
{1+\frac{1}{m}}{Z}
\end{equation}
Let $f_0(Z)$ be 
the quotient of these two functions (a ``Schwarz map" for $H_0$).
We evaluate $f_0(\phi_0(z))$ using suitable hypergeometric
identities from the data base in Section \ref{database}.
Up to a constant multiple, the Schwarz maps are equal to $z^{-1/m}$.
\algstep {\bf Step 4.} Consider two hypergeometric solutions of $H_1$,
\begin{equation}
G_1=\hpg21{a,b\,}{c}{X},\qquad 
G_2=X^{1-c}\,\hpg21{1+a-c,1+b-c}{2-c}{X},
\end{equation}
corresponding to the functions in (\ref{twoh1sols}). The local exponent differences $1-c$,
$c-a-b$ and (say) $a-b$ are the positive local exponent differences of the given equation.
In particular,  $c<1$ and $c-a-b>0$. Besides, $c$ is a rational number with the denominator $m$. 
The functions $F_1$ and $G_1$ are identified under Klein's pull-back transformation
up to a power factor. The functions $F_2$ and $G_2$ are identified up to the same power factor
and a constant multiple.
\algstep  {\bf Step 5.} In the data base of Section \ref{database} there are four
evaluations of hypergeometric functions with the Schwarz type of $H_1$. 
Possibly up to Euler's transformation \cite[(2.2.7)]{specfaar} 
two of them are contiguous to $G_1$, and the other two are contiguous to $G_2$. 
We find a contiguous relations between $G_1$ (or its Euler's transformation) 
and the former two data base functions, and a contiguous relation between $G_2$ 
(or its Euler's transformation) and the latter two data base functions. 
We evaluate the functions using the contiguous relations and the data base evaluations.
The data base evaluations introduce radical expressions in $x,\xi$. 
Coefficients of the contiguous relations are rational functions of $X$; 
they should not be evaluated in terms of $x,\xi$. The resulting expressions of 
the functions $G_1$ and $G_2$ have the following form: a fractional power expression in $x$
and possibly $\xi$ (which can be standardized for each Schwarz type),
times a rational function in $X,x$ and possibly $\xi$ (of bounded degree in $x,\xi$).
\algstep {\bf Step 6.} Let $f_1(X)$ denote the Schwarz map $G_1/G_2$ for $H_1$.
We consider $f_1(X)$ explicitly evaluated like the functions $G_1,G_2$ in the previous step:
a standardized fractional power function in $x,\xi$ and a rational function in $X,x,\xi$.
The Schwarz maps  $f_0(Z)$ and $f_1(X)$ are identified up to a constant multiple $w$
by Klein's pull-back covering. The identification gives a rational relation $z=w\,\Phi(X,x,\xi)$.
To find the constant $w$, we evaluate this identity at the following pairs of $z$- and $(x,\xi)$-points,
lying above $Z=0$ and $X=0$, respectively. 
In the icosahedral cases, we have two possibilities for the $z$-point; the right possibility
is determined by checking compatibility of the given automorphisms of the $z$- and $(x,\xi)$
Darboux coverings.
\begin{itemize}
\item For the tetrahedral type $(1/2,1/3,1/3)$: evaluation at $x=-4$, or $z=-4$.
\item For the tetrahedral type $(1/3,1/3,2/3)$: evaluation at $x=-2$.
\item For the octahedral type $(1/2,1/3,1/4)$: evaluation at $x=1$, or $z=1$.
\item For the octahedral type $(2/3,1/4,1/4)$: evaluation at $x=-1$.
\item For the icosahedral types $(1/2,1/3,1/5)$ and $(1/2,1/3,2/5)$:\\
evaluation at $x=(11+5\sqrt{5})/2$, or $z=(11\pm 5\sqrt{5})/2$;\\ 
the automorphism $x\mapsto -1/x$, or $z\mapsto -1/z$.
\item For the icosahedral type $(1/2,1/5,2/5)$:
evaluation at $x=(1+\sqrt{5})/2$;\\
the automorphism $x\mapsto -1/x$.
\item For the icosahedral types $(1/3,1/3,2/5)$ and $(1/3,2/3,1/5)$:\\ 
evaluation at $(x,\xi)=\left((11+5\sqrt{5})/6,0\right)$;
the automorphism $x\mapsto -1/9x$,  $\xi\mapsto \xi/9x^2$.
\item For the icosahedral types $(2/3,1/5,1/5)$ and $(1/3,2/5,3/5)$:\\ 
evaluation at $(x,\xi)=\left((-3+7/\sqrt{5})/2,7-3\sqrt5\right)$;\\
the automorphism $x\mapsto -1/5x$, $\xi\mapsto \xi/5x^2$.
\item For the icosahedral type $(1/3,1/5,3/5)$: \\
evaluation at $(x,\xi)=\left((-3+\sqrt{5})/8,(-5+3\sqrt{5})/8\right)$;\\
the automorphism $x\mapsto 1/16x$, $\xi\mapsto -\xi/16x^2$.
\item For the icosahedral types $(1/5,1/5,4/5)$ and $(2/5,2/5,2/5)$:\\ 
evaluation at $(x,\xi)=\left((1+\sqrt{5})/2,0\right)$;
the automorphism $x\mapsto -1/x$, $\xi\mapsto \xi/x^2$.
\end{itemize}
{\bf Step 7.} We have three algebraic equations produced by Step 2, and the equation
$z=w\,\Phi(X,x,\xi)$ produced by Ste p 6.  The variables are $Z$, $X$, $z$, $x$ and $\xi$. 
The equations have bounded degree in $z$, $x$, $\xi$. Using Gr\"obner bases techniques (or
resultants) we eliminate $z$, $x$ and (if relevant) $\xi$, and get a polynomial relation between
$X$ and $Z$. Klein's theorem ensures this polynomial relation has a factor linear in $Z$,
giving the solution $Z=\Phi(X)$. \qed
\end{algorithm}

This algorithm is demonstrated by examples in Sections \ref{sec:examples} and 
Section \ref{sec:ixamples}. Here we give several comments to some steps.

Let $(e_0,e_1,e_\infty)$ denote the positive local exponent differences of the given equation $H_1$.
The degree $d$ of Klein's pull-back covering can be computed using part 2 of 
\cite[Lemma 2.5]{algtgauss}:
\begin{equation}
d=\frac{e_0+e_1+e_\infty-1}{\frac12+\frac13+\frac1m-1}
= \frac{6m}{6-m}(e_0+e_1+e_\infty-1).
\end{equation}
This formula has an interpretation as the area ratio of the Schwarz triangles for the given and standard
hypergeometric equations. 

Steps 0 and 3 are the same for input differential equations with same mondromy group.
Steps 1, 2 and probably 6 can be set up to work in the same way for the input of the same Schwarz type.

In Step 3, let $s=z^{-1/m}$ denote the Schwarz maps for a standard tetrahedral, octahedral and icosahedral equation. The inverse Schwarz map is a rational function $S_m(s)$ of $s$ 
(with $m=3,4,5$). We recognize classical invariants of the Galois groups in these expressions:
\begin{eqnarray*}
&& S_3(s)=\frac{s^3\left(s^3+4\right)^3}{4\left(2s^3-1\right)},\qquad
S_4(s)=\frac{108\,s^4\left(s^4-1\right)^4}{\left(s^8+14s^4+1\right)^3},\\
&& S_5(s)=\frac{1728\,s^5\left(s^{10}-11s^5-1\right)^5}
{\left(s^{20}+228s^{15}+494s^{10}-228s^5+1\right)^3}.
\end{eqnarray*}

In Step 4, the parameters $a,b,c$ have the following expressions in terms of 
the local exponent differences:
\begin{equation}
a=\frac{1-e_0-e_1-e_\infty}2,\qquad a=\frac{1-e_0-e_1+e_\infty}2, \qquad
c=1-e_0.
\end{equation}

Euler's fractional-linear transformation \cite[(2.2.7)]{specfaar} in Step 5 reads
\begin{equation} \label{eq:euler}
\hpg21{a,b}{c}{\,z\,} = (1-z)^{c-a-b}\,\hpg21{c-a,\,c-b}{c}{\,z\,}.
\end{equation}
The function $G_1$ (or its Euler's transformation) might be contiguous to the last two
of the four evaluations of the relevant Schwarz type in the data base,
while the function $G_2$ (or its Euler's transofrmation) would be contiguous
the the first two data base evaluations; this is demonstrated in Example \ref{ex:tetr2}.
There are (at most) 4 different contiguity orbits for $G_1$ (and $G_2$) for each Schwarz type;
for each of those orbits, the constant $w$ in the next step can probably be determined
by a fixed evaluation of the relation $z=w\Psi(X,x,\xi)$ even in icosahedral cases.

Note that each pair of contiguous evaluations in the data base is presented with
the same fractional-power part, in the cases of genus 1 Darboux curves as well.
This allows the algorithm to keep the same fractional-power part in all contiguous evaluations,
keeping the same pattern of $(x,\xi)$ variables for Step 7 in each Schwarz case.

In Step 6, we evaluate the relation $z=w\Psi(X,x,\xi)$ at a $(x,\xi)$-point above $X=0$ (but not $x=0$)
and at a $z$-point above $Z=0$ (but not $z=0$). Analysis of diagram (\ref{di:darboux}) for each 
Schwarz type shows that the $(x,\xi)$ point lies above the $z$-point, except that we may have a
choice of two $z$-values in the icosahedral cases. We check the two $z$-values by the given
automorhisms of Darboux curves. On the other hand, just compatibility of the given automorphisms 
determines the constant multiple $w$ up to a sign. Evaluation at $X=1$, $Z=1$ using Gauss formula
\cite[Theorem 2.2.2]{specfaar} can be used, but the $X=1$ value may need to be twisted to
a correct branch (which probably depends only on the contiguity orbit of $G_1$);
see Example \ref{ex:tetr1}.

In Step 7, particular elimination tricks (while using a modern computer algebra system)
are essentially needed only in icosahedral cases.
Examples suggest that elimination using Gr\"obner bases or resultants returns 
an integer power of $Z-\Phi(X)$, with cleared denominator.
The icosahedral relations have degree 12 (or a little higher) in the variable $x$,
which is hard to eliminate directly, especially together with $\xi\neq 0$ on genus 1
Darboux curves. The respective automorphism symmetry of Step 6 can be used
to half the degree in $x$. The new variable would look like $\alpha=x-1/x$ or so; 
see Example \ref{ex:icosa1}. As the relation $z=w\Psi(X,x,\xi)$ and the output have highly 
factorizable forms, elimination can be done fast for each factor, and the resulting factors in $X$
can be combined and then normalized using evaluation at $X=1$ or differentiation (to find
ramification points that must be above $Z=1$); see Example \ref{ex:icosa1}. 
Algebraic computations can probably be streamlined for each Schwarz type 
using a normal series of the solvable monodromy group and respective resolvents, 
with simplification at each step to a variable of higher invariance. 

\section{Tetrahedral examples}
\label{sec:examples}

Here we demonstrate Algorithm \ref{kleinalg} on several hypergeometric
equations with the tetrahedral monodromy group. 
The standard hypergeometric equation for this monodromy group 
has the local exponent differences $(1/2,1/3/,1/3)$.
Following Step 0, we assign the local exponent difference $1/3$ 
to the points $Z=0$ and $Z=\infty$. The local exponent difference $1/2$
is assigned to the point $Z=1$. The Darboux covering $Z=\phi_0(z)$  in Step 2 is
\begin{equation} \label{darbtrsd}
Z=\frac{z\,(z+4)^3}{4(2z-1)^3}. 
\end{equation}
The Schwarz map in Step 3 is  
\begin{equation} \label{schwartzd}
f_0(Z)=\hpg{2}{1}{\frac14,\,-\frac1{12}\,}{\frac23}{Z}\left/
Z^{1/3}\,\hpg{2}{1}{\frac14,\;\frac7{12}\,}{\frac43}{Z}\right. = \left(
-\frac{1}{16\,z} \right)^{1/3},
\end{equation}
evaluated using  (\ref{darbtrsd}) and (\ref{fptetra1}), (\ref{fptetra1b}).

The first three examples are of the same Schwarz type $(1/2,1/3,1/3)$. 
In accordance with Step 1, the points $X=0$, $X=1$, $X=\infty$ are assumed to lie
above the points $Z=0$, $Z=1$, $Z=\infty$, respectively.
The Darboux curve for those hypergeometric equations is rational, hence $\xi=0$.
The Darboux covering $X=\phi_1(x)$ in Step 2 mimics (\ref{darbtrsd}):
\begin{equation} \label{darbtrst3}
X=\frac{z\,(z+4)^3}{4(2z-1)^3}. 
\end{equation}

\begin{example} \rm \label{ex:tetr1}
We compute Klein's pull-back covering for hypergeometric equations
with the local exponent differences $(1/2,1/3,2/3)$. The degree is expected to be 3.

The Schwartz map in Step 4 is
\begin{equation} \label{eq:schwt3a}
\left.\hpg21{-\frac14,\,\frac{5}{12}}{\frac23}{X} \right/ X^{1/3}\hpg21{\,\frac34,\;\frac{1}{12}}{\frac43}{X},
\end{equation}
The two $\hpgo21$ functions are not contiguous to the data base functions
(\ref{fptetra1})--(\ref{fptetra1z}), but by Euler's transformation (\ref{eq:euler}) 
they are equal, respectively, to
\[
\sqrt{1-X}\,\hpg21{\frac14,\,\frac{11}{12}}{\frac23}{X}, \qquad
\sqrt{1-X}\,\hpg21{\frac54,\,\frac7{12}}{\frac43}{X}. 
\]
We have the following contiguous relations for these functions,
with $(a,b,c)=(1/4,-1/12,2/3)$ or $(1/4,7/12,4/3)$, respectively:
\begin{eqnarray}
\hspace{-6pt} bc\,(1\!-\!X)\,\hpg21{\!a,b\!+\!1}{c}{X\!} \!\equal
\!c\left((a\!-\!b)X\!+\!b\right)\hpg21{\!a,b}{c}{X\!}\!+a(b\!-\!c)X\,\hpg21{\!a\!+\!1,b}{c\!+\!1}{X\!}\!,
\hspace{8pt}\\ \label{eq:contmm}
a\,(1\!-\!X)\,\hpg21{\!a\!+\!1,b}{c}{X\!} \!\equal
(a+1-c)\,\hpg21{\!a,b}{c}{X}+(c-1)\,\hpg21{\!a,b\!-\!1}{c\!-\!1}{X}.
\end{eqnarray}
Using data base formulas (\ref{fptetra1})--(\ref{fptetra1z}), Step 5 returns
\begin{eqnarray*}
\hpg21{-\frac14,\,\frac{5}{12}}{\frac23}{\frac{x\,(x+4)^3}{4(2x-1)^3}}\equal (1+x)\,(1-2x)^{-3/4}, \\
\hpg21{\,\frac34,\;\frac{1}{12}}{\frac43}{\frac{x\,(x+4)^3}{4(2x-1)^3}}\equal
\frac1{1+\frac14x}\,(1-2x)^{1/4}.
\end{eqnarray*}
In terms of $x$, the Schwarz map $f_1(X)$ in (\ref{eq:schwt3a}) is a constant times $x^{-1/3}(x+1)$.
This must be proportional to (\ref{schwartzd}).  Hence $z=w\,x/(x+1)^3$, 
where the constant $w=-27$ is obtained by evaluation $(x,y)=(-4,-4)$. 
In the final step, we eliminate $x,z$ using (\ref{darbtrst3}) and (\ref{darbtrsd})
by whatever straightforward method, and obtain the result $Z=27X/(4X-1)^3$.
One can try and note that computations become ugly when $w$ is taken wrongly.

The cubic transformation $Z=27X/(4X-1)^3$ is the same as in a classical cubic transformation
of Gauss hypergeometric functions.

Direct comparison of Schwarz maps (\ref{schwartzd}) and (\ref{eq:schwt3}) by their evaluations at,
respectively, $Z=1$ and $X=1$, respectively, does not a right identification. Evidently,
as $X$ varies from 0 to 1, $Z$  leaves the convergence disk, goes through the infinity,
and approaches $1$ on other branch. The explicit evaluations at $Z=1$ and $X=1$ can be done
using Gauss formula and the data in \cite{gammaeval}. The values are, respectively,
\begin{eqnarray*}
\frac{\Gamma(2/3)\,\Gamma(13/12)}{\Gamma(4/3)\,\Gamma(5/12)} =
\frac{\sqrt{3}+1}4,\qquad
\frac{\Gamma(2/3)\,\Gamma(7/12)}{4\,\Gamma(4/3)\,\Gamma(11/12)} =
\frac{3\,(\sqrt{3}-1)}4.
\end{eqnarray*}
The right identification is achieved after conjugating $\sqrt3$; this twist is probably preserved
under contiguous relations and can be analyzed for each Schwartz type.
\end{example}

\begin{example} \rm \label{ex:tetr2}
Here we compute a degree 5 Klein's pull-back covering for the local exponent differences 
$(1/2,2/3,2/3)$.  Up to a constant multiple, the Schwartz map in Step 4 is
\begin{equation} \label{eq:schwt3b}
\left.\hpg21{\frac14,\,-\frac{5}{12}}{\frac13}{X} \right/ 
X^{2/3}\hpg21{\,\frac14,\;\frac{11}{12}}{\frac53}{X},
\end{equation}
The two functions are contiguous to the data base evaluations (\ref{fptetra1})--(\ref{fptetra1z}),
though the first function is in the contiguous orbit of the later two evaluations
(\ref{fptetra1b})--(\ref{fptetra1z}), and vice versa. The numerator function is actually presented in
(\ref{fptetra1z}). For the denominator $\hpgo21$ function we have this contiguous relation
with $(a,b,c)=(1/4,-1/12,2/3)$:
\begin{eqnarray*}
b\,(a-c)\,\hpg21{\!a,b\!+\!1}{c+1}{X} \equal
c\,(a-b)\,\hpg21{\!a,b}{c}{X}+a\,(b-c)\,\hpg21{\!a\!+\!1,b}{c\!+\!1}{X}.
\end{eqnarray*}
Using data base formulas (\ref{fptetra1})--(\ref{fptetra1a}), Step 5 returns
\begin{eqnarray*}
\hpg21{\frac14,\,\frac{11}{12}}{\frac53}{X} \equal
\frac{1-x/5}{1+x/4}\,(1-2x)^{3/4}.
\end{eqnarray*}
In terms of $x$, the Schwarz map $f_1(X)$ in (\ref{eq:schwt3b}) is a constant times 
$x^{-2/3}(1+5x/2)/(1-x/5)$.  Therefore $z=w\,x^2(x-5)^3/(5x+2)^3$, where the constant 
$w=-2$ is obtained  by evaluation at $(x,z)=(-4,-4)$. 
The final step returns $Z=-X^2(4X-5)^3/(5X-4)^3$; a modern computer algebra system
has little problem eliminating $x,z$.

The following transformation of Gauss hypergeometric functions follows: 
\begin{eqnarray} \label{eq:hpg5}
\hpg21{\frac14,\,-\frac{5}{12}}{\frac13}{X}= { \left(1-\frac{5X}4 \right)^{1/4}}\;
\hpg{2}{1}{\frac14,\,-\frac1{12}\,}{\frac23}{-\frac{X^2(4X-5)^3}{(5X-4)^3}}. \\ \nonumber
\end{eqnarray}
\end{example}

\begin{example} \rm
Now we compute a degree 7 Klein's pull-back covering for the local exponent differences 
$(3/2,1/3,1/3)$. Up to a constant multiple, the Schwartz map in Step 4 is
\begin{equation} \label{eq:schwt3}
\left.\hpg21{-\frac14,\,-\frac{7}{12}}{\frac23}{X} \right/ 
X^{1/3}\hpg21{-\frac14,\;\frac{1}{12}}{\frac43}{X},
\end{equation}
Like in Example \ref{ex:tetr1}, we have contiguous relations for Euler transformations 
of these functions. The contiguous relations have to be applied,
with $(a,b,c)=(1/4,-1/12,2/3)$ or $(1/4,7/12,4/3)$, respectively:
\begin{eqnarray*}
b\,c\,(1\!-\!X)^2\,\hpg21{\!a\!+\!1,b\!+\!1}{c}{X} \equal
c\,\left((a-c+1)X+b\right)\,\hpg21{\!a,b}{c}{X}\\
&&+(a+b+1-c)(b-c)X\,\hpg21{\!a+1,b}{c+1}{X},\\
a\,b\,(1\!-\!X)^2\,\hpg21{\!a\!+\!1,b\!+\!1}{c}{X} \equal
(a-c+1)\,\left(aX+b-c+1\right)\,\hpg21{\!a,b}{c}{X}\\
&&+(a+b+1-c)(c-1)\,\hpg21{\!a,b\!-\!1}{c\!-\!1}{X}.
\end{eqnarray*}
Using data base formulas (\ref{fptetra1})--(\ref{fptetra1z}), Step 5 returns
\begin{eqnarray*}
\hpg21{-\frac14,\,-\frac{7}{12}}{\frac23}{\frac{x\,(x+4)^3}{4(2x-1)^3}}\equal
{\textstyle\left(1-7x-\frac72x^2\right)}\,(1-2x)^{-7/4}, \\
\hpg21{-\frac14,\;\frac{1}{12}}{\frac43}{\frac{x\,(x+4)^3}{4(2x-1)^3}}\equal
\frac{1-x-\frac{1}{14}x^2}{1+\frac14x}(1-2x)^{-3/4}.
\end{eqnarray*}
In Step 6, we conclude that $z=w\,x\,(x^2+14x-14)^3/(7x^2+14x-2)^3$, 
where the constant $w=-1$ by evaluation $(x,y)=(-4,-4)$.  The final step gives 
\begin{equation}
 Z = -\frac{X(X^2-42X-7)^3}{(7X^2+42X-1)^3}.
 \end{equation}
Similarly to (\ref{eq:hpg5}), 
the following transformation of Gauss hypergeometric functions follows: 
\begin{eqnarray}
\hpg21{-\frac14,\,-\frac{7}{12}}{\frac23}{X}= { \left(1-42X-7X^2 \right)^{1/4}}\;
\hpg{2}{1}{\frac14,\,-\frac1{12}\,}{\frac23}{-\frac{X(X^2-42X-7)^3}{(7X^2+42X-1)^3}}.\\ \nonumber
\end{eqnarray}
\end{example}

Here are a few more Klein's pull-back coverings for hypergeometric equations
of Schwarz type $(1/2,1/3,1/3)$, of degree 7 or 9. The given are local exponent differences
and Klein's pull-back coverings:
\begin{eqnarray*}
(1/2,1/3,4/3): && Z = -\frac{X(256X^2-448X+189)^3}{27(28X-27)^3},\\
(1/2,2/3,4/3): && Z = \frac{19683X^2(4X-1)^3}{(256X^3-192X^2+21X-4)^3},\\
(1/2,1/3,5/3): && Z = -\frac{19683X(128X-125)^3}{(16384X^3-30720X^2+14880X-625)^3},\\
(3/2,1/3,2/3): && Z = -\frac{729X(5X^2+14X+125)^3}{(4X^3+15X^2-690X-625)^3}.
\end{eqnarray*}

Next we give Klein's pull-back covering of degree 14 for a hypergeometric equation of 
the other tetrahedral Schwarz type $(1/3,1/3,2/3)$,
underscoring the strategy of Step 5 to keep coefficients of contiguous relations
as functions of $X$ only, so to avoid higher degree in the variable $x$ to be eliminated.
The simplest Klein's pull-back covering for Schwarz type $(1/3,1/3,2/3)$ is a quadratic transformation.
Klein's coverings for  
the local exponent differences $(2/3,2/3,2/3)$ or $(1/3,1/3,4/3)$ are 
specializations of classical transformations of degree 6 for Gauss hypergeometric functions. 
Here are two Klein's coverings of degree 10:
\begin{eqnarray*}
(1/3,2/3,5/3): && Z = \frac{4X(256X^3-640X^2+520X-135)^3}{27(X-1)^2(32X-27)^3},\qquad \\
(2/3,2/3,4/3): && Z = -\frac{X^2(X-1)^2(16X^2-16X+5)^3}{4(5X^2-5X+1)^3}.
\end{eqnarray*}

\begin{example} \label{exatetra} \rm
Here we compute Klein's pull-back covering for the local exponent differences 
$(2/3,4/3,4/3)$. We keep the same setting of Step 0, but we have to adjust Step 1.
Eventually, the ramification pattern must be $3+3+3+3+2=4+4+3+3=2+2+2+2+2+2+2$,
using the notation of \cite{algtgauss}. We assign the local exponent $4/3$ 
to the points $X=0$ and $X=1$, and the local exponent $2/3$ to the point $X=\infty$.
Since the point $X=0$ lies above $Z=0$ by assumption, the point $X=1$ will end up there as well.

In Step 2 the Darboux covering for the standard equation is the same (\ref{darbtrsd}), 
while the Darboux covering for the given equation $H_1$ is 
\begin{equation} \label{darbtrst4}
X=\frac{x\,(x+2)^3}{(2x+1)^3}.
\end{equation}
Similarly, in Step 3 the Schwarz map $f_0(Z)$ for the standard equation is the same (\ref{schwartzd}),
while the Schwarz map $f_1(X)$ in Step 4 is
\begin{equation} \label{schwartzx}
\hpg{2}{1}{-\frac12,\,-\frac76}{-\frac13}{X}\left/
X^{4/3}\,\hpg{2}{1}{\frac16,\,\frac56\,}{\frac73}{X}. \right.
\end{equation}
Contiguous relations for these functions are the following, 
with $(a,b,c)=(1/2,-1/6,2/3)$ or $(1/6,5/6,4/3)$, respectively. 
\begin{eqnarray*}
(1-c)\,\hpg21{\!a\!-\!1,b\!-\!1}{c\!-\!1}{X\!} \!\equal
\!\left((a\!-\!1)X\!+\!b\!-\!c\!+\!1\right)\hpg21{\!a,b}{c}{X\!}+b(1\!-\!X)\hpg21{\!a,b\!+\!1}{c}{X\!},\\
(c\!-\!a)(c\!-\!b)X\hpg21{\!a,\,b}{c\!+\!1}{X\!} \!\equal
c\left((c\!-\!b)X\!-\!c\!+\!1\right))\,\hpg21{\!a,b}{c}{X}+c\,(c\!-\!1)\,\hpg21{\!a,b\!-\!1}{c\!-\!1}{X},
\end{eqnarray*}
Utilizing data base formulas (\ref{fptetra2})--(\ref{fptetra2z}) we get the expressions in Step 5:
\begin{eqnarray} \label{fptetra2p}
\hpg{2}{1}{\!-\frac12,-\frac76\,}{-\frac13}{X} \equal \left( \frac{1-3X}{2}+
\frac{1-X}{2}\,\frac{(1+2x)^2}{(1-x)^2}\right)\frac{1}{\big(1+2x\big)^{1/2}}.\\
\label{fptetra2q}\hpg{2}{1}{\frac16,\;\frac56\,}{\frac73}{X}  \equal 
\frac{8}{21} \left( \frac{3X-2}{X}
+\frac{2}{X}\,\frac{1+\frac12x}{1+2x}\right)
\frac{\big(1+2x)^{1/2}\,\big(1+x\big)^{1/3}}{1+\frac12x}.
\end{eqnarray}
In Step 6 we compute Schwarz map (\ref{schwartzx}):
\begin{equation}
f_1(X,x)=\frac{21\,(x+2)}{32\,X^{1/3}\,(1+x)^{1/3}(1-x)^2}\,
\frac{(1-3X)(1-x)^2+(1-X)(1+2x)^2}{(3X-2)(1+2x)+(x+2)}.
\end{equation}
Further evaluation using (\ref{darbtrst4}) would probably give a rather nice expression in $x$ only,
but we wish to avoid higher degree in $x$. We may only wish to simplify the power $X^{1/3}$
in the denominator. The Schwarz map $f_1(X,x)$ has to be identified with $f_1(Z)$;
we conclude that
\begin{equation}\label{schwart3w}
\frac{1}{z}=\frac{w\;(2x+1)^{3}}{x\,(1+x)(1-x)^6}\, \left(
\frac{(1-3X)(1-x)^2+(1-X)(1+2x)^2}{(3X-2)(1+2x)+(x+2)} \right)^3
\end{equation}
for some constant $w$. To find the constant, we evaluate at $z=-4$, $x=-2$ and $X=0$;
this is compatible with (\ref{darbtrst4}), and there are no zeroes in multiplicative terms.
We find $w=1/2$.

In Step 7 we can use straightforward means (within a modern computer algebra system like {\sf Maple})
to eliminate $x$ and $z$, and get the following unique result:
\begin{equation}
Z=-\frac{108\;X^4\;(X-1)^4\;(27X^2-27X+7)^3}
{(189X^4-378X^3+301X^2-112X+16)^3}.
\end{equation}
\end{example}

\section{Two icosahedral examples}
\label{sec:ixamples}

Computation of Klein's pull-back coverings for icosahedral hypergeometric equations
is a level harder than for tetrahedral or octahedral hypergeometric equations. As Darboux
coverings have degree 12 and Darboux curves can be non-rational, straightforward elimination
medthods are inadequate. The following two examples demonstrate some of simplification 
techniques suggested at the end of Section \ref{kleinmorphs}.

The standard hypergeometric equation for the icosahedral monodromy group 
has the local exponent differences $(1/2,1/3/,1/5)$. We assign these local exponents to, 
respectively, $Z=1$, $Z=\infty$ and $Z=0$.
The Darboux covering $Z=\phi_0(z)$ in Step 2 is given in (\ref{isophi1}):
\begin{equation} \label{darbtrsdi}
Z=\frac{1728\;z\;(z^2-11z-1)^5}{(z^4+228z^3+494z^2-228z+1)^3}.
\end{equation}
The Schwarz map in Step 3 is  
\begin{equation} \label{schwartzdi}
f_0(Z)=\hpg{2}{1}{\frac{19}{60},\,-\frac1{60}\,}{\frac45}{Z}\left/
Z^{1/5}\,\hpg{2}{1}{\frac{11}{60},\;\frac{59}{60}\,}{\frac65}{Z}\right. = \left(
-\frac{1}{1728\,z} \right)^{1/5},
\end{equation}
evaluated using  (\ref{darbtrsdi}) and (\ref{fpicosa1}), (\ref{fpicosa1b}).

\begin{example} \rm \label{ex:icosa1}
Here we  compute a degree 11 pull-back covering for hypergeometric equations
with the local exponent differences $(1/2,2/3,1/5)$. In Step 1, we choose
 the points $X=0$, $X=1$, $X=\infty$ to lie above the points $Z=0$, $Z=1$, $Z=\infty$, respectively.
The Darboux covering $X=\phi_1(x)$ in Step 2 mimics (\ref{darbtrsdi}):
\begin{equation} \label{darbtrsdix}
X=\frac{1728\;x\;(x^2-11x-1)^5}{(x^4+228x^3+494x^2-228x+1)^3}.
\end{equation}
The Schwartz map in Step 4 is
\begin{equation} \label{eq:schwi11}
\left.\hpg21{\frac{29}{60},\,-\frac{11}{60}}{\frac45}{X} \right/ 
X^{1/5}\hpg21{\frac1{60},\;\frac{41}{60}}{\frac65}{X},
\end{equation}
An explicit expression for Euler's transformation of the numerator function is presented directly
in (\ref{fpicosa1a}), while Euler's transformation for the denominator function can be computed
from contiguous relation (\ref{eq:contmm}) with $(a,b,c)=(11/60,31/60,6/5)$ and data base formulas 
(\ref{fpicosa1b})--(\ref{fpicosa1z}).  Step 5 returns the explicit expressions
\begin{eqnarray*}
\hpg21{-\frac{11}{60},\,\frac{29}{60}}{\frac45}{X}\equal
(1+66x-11x^2)\,(1-228x+494x^2+228x^3+x^4)^{-11/20}, \\
\hpg21{\frac{1}{60},\;\frac{41}{60}}{\frac65}{X}\equal
\frac{1+6x-\frac{1}{11}x^2}{1+11x-x^2}\,(1-228x+494x^2+228x^3+x^4)^{1/20}.
\end{eqnarray*}
Up to a constant multiple, we have $z=w\,x(11+66x-x^2)^5/(1+66x-11x^2)^5$ in Step 6. 
The constant $w=-1$ is determined by the evaluation $x=z=(11+5\sqrt{5})/2$.
The alternative evaluation $z=(11-5\sqrt{5})/2$ gives $w=(123-55\sqrt5)/2$,
but then the relation between $x$ and $z$ is not invariant under the given automorphisms 
$x=-1/x$, $z\mapsto -1/z$.

To simplify equations (\ref{darbtrsdi}), (\ref{darbtrsdix}) and the relation for $z$ with $w=-1$,
we introduce the variables $\alpha=x-1/x$ and $\beta=z-1/z$. Due to the symmetry between
higher and lower coefficients in the polynomial factors in (\ref{darbtrsdi}) and (\ref{darbtrsdix}),
we can ``replace" $x$, $z$ by the variables $\alpha,\beta$; we get the simpler expressions
\begin{equation} \label{eq:xyab}
X=\frac{1728(\alpha-11)^5}{(\alpha^2+228\alpha+496)^3}, \qquad
Z=\frac{1728(\beta-11)^5}{(\beta^2+228\beta+496)^3}.
\end{equation}
The expression for $z$ with $w=-1$ is replaced by
\[
\beta=z-\frac1{z}=-\frac{\alpha^{11}-660\alpha^{10}-25937228680\alpha^9+\ldots}
{(11\alpha^2-792\alpha+4256)^5}.
\]
Elimination of $\alpha$, $\beta$ is rather straightforward now. The result is
\begin{eqnarray}
Z=\frac{X\,(102400X^2-11264X-11)^5}{(180224000X^3+4325376X^2-21252X+1)^3}.\\ \nonumber
\end{eqnarray}
\end{example}

\begin{example} \rm \label{ex:icosa2}
Here we compute a degree 18 pull-back covering for hypergeometric equations
with the local exponent differences $(1/5,1/5,6/5)$. The Darboux curve is of genus 1,
defined in (\ref{darbouxc4}). In Step 1, we copy the Darboux covering from (\ref{isophi6}):
\begin{equation} \label{eq:darb18}
X=\frac{16\;\xi\,\left(1+x-x^2\right)^2\,(1-\xi)^2}{(1+\xi+2x)\,(1+\xi-2x)^5}, 
\qquad \xi^2=x\,(1+x-x^2).
\end{equation}
In Step 3, the Schwartz map $f_1(X)$ is
\begin{equation} \label{eq:schwi18}
\left.\hpg21{\frac{9}{10},\,-\frac{3}{10}}{\frac45}{X} \right/ 
X^{1/5}\hpg21{\frac{11}{10},-\frac{1}{10}}{\frac65}{X},
\end{equation}
Contiguous relations for these functions are the following, 
with $(a,b,c)=(7/10,-1/10,4/5)$ or $(1/10,9/10,6/5)$, respectively:
\begin{eqnarray*}
c\;\hpg21{\!a\!-\!1,b\!+\!1}{c}{X} \!\equal
c\;\hpg21{\!a,b}{c}{X}+(a\!-\!b\!-\!1)X\,\hpg21{\!a,b+1}{c+1}{X},\\ 
a(b\!-\!c)\hpg21{\!a\!+\!1,b\!-\!1}{c}{X\!} \!\equal
\!(b\!-\!1)(a\!\!-\!c+\!1)\,\hpg21{\!a,b}{c}{X\!}+(c\!-\!1)(b\!-\!a\!-\!1)\hpg21{\!a,b\!-\!1}{c\!-\!1}{X\!}.
\end{eqnarray*}
Utilizing data base formulas (\ref{icosellipta6})--(\ref{icosellipta6z}) we get these expressions
in Step 5:
\begin{eqnarray} 
\hpg{2}{1}{\frac9{10},-\frac3{10}}{\frac45}{X} \equal \frac{
(1-2x-3y)\,(1-\xi+2x)^{1/15}\,}{(1-\xi)^{7/5}(1+\xi+2x)^{7/30}\,(1+\xi-2x)^{3/2}},\\
\hpg{2}{1}{\frac{11}{10},\,-\frac1{10}}{\frac65}{X} \equal \frac{(3\xi-2x-x^2)
\;(1+\xi)^{1/10}\,(1-\xi)^{3/10}}{3\xi\,(1-\xi+2x)^{1/30}
(1+\xi+2x)^{2/15}\sqrt{1+\xi-2x}}.
\end{eqnarray}
Step 6 gives a Schwarz map $f_1(x,\xi)$ and the following relation: 
\begin{equation}
z=\frac{w\,(1-x)(3\xi-2x-x^2)^5}{x^2(1+x)(3\xi+2x-1)^5},
\end{equation}
where $w=1$ by evaluation at $(x,\xi)=\left((1+\sqrt5)/2, 0\right)$ and $z=(11+5\sqrt{5})/2$.
We use the compatible automorphsims $(x,\xi)\mapsto (-1/x,\xi/x^2)$ and $z\mapsto-1/z$ of Darboux coverings to check the value of $w$, and simplify equations for Step 7. Correspondingly, we introduce the new variables
\[
\alpha=x-\frac1x,\qquad \zeta=\frac{\xi\,(x^2+1)}{x^2}, \qquad \beta=z-\frac1z.
\]
Equations (\ref{eq:darb18}) can be simplified to
\[
\frac1X = \frac12-\frac{\zeta\,(\alpha^2-22\alpha-4)}{16(\alpha-1)^3},
\qquad \zeta^2=(1-\alpha)(4+\alpha^2).
\]
From here we can eliminate $\zeta$ linearly and get the following equation of degree 6 for $\alpha$:
\begin{equation}
\alpha(\alpha+4)^5X^2-256(\alpha-1)^5X+256(\alpha-1)^5=0
\end{equation}
Alternatively, one can express all cubic monomials in $x,\xi$ in lower degree monomials
by a Gr\"obner basis. That would be a basis to express
\[
\beta=z-\frac1z=\frac{-6\,\zeta\,(5\alpha^7+65150\alpha^6+\ldots)+(\alpha^9+58662\alpha^8+\ldots)}
{\alpha\,(3\zeta+11\alpha-6)^5}
\]
and a few more (a priori specified) 
increasingly invariant expressions towards $Z$; we would keep
at most 6 low degree monomials in $\alpha,\zeta$. 
But automatic computation from this point with Gr\"obner basis in all variables does not look feasible.
Without a resolvent sequence, a practical trick is to take (the numerators of) the factors $\beta-11$
and $\beta^2+228\beta+496$ in (\ref{eq:xyab}), ``eliminate" $\alpha,\zeta$ from them using resultants, 
and recognize suitable polynomial factors in $X$ for a possible final expression, keeping in mind its
whole ramification pattern.
Candidate expressions can be normalized by the condition that other ramification points must lie above 
$Z=1$. In this case, we get a unique candidate factor $512X^2-512X+3$ for the numerator, and
a unique factor in $X$ for the denominator, and an additional ramification point is $X=1/2$. 
The final result is
\begin{equation}
Z=-\frac{108\,X\,(X-1)\,(512X^2-512X+3)^5}
{(1048576X^6\!-\!3145728X^5\!+\!3244032X^4\!-\!1245184X^3\!+\!94848X^2\!+\!3456X\!+\!1)^3}.
\end{equation}
\end{example}

\section{A data basis of hypergeometric identities}
\label{database}

Here we present a data base of Darboux coverings and Darboux evaluations
of algebraic hypergeometric solutions with small local exponent
differences. For each Schwarz type we evaluate four functions:
two hypergeometric solutions of the main representative hypergeometric equation, 
and a contiguous ``companion" function to each of these two functions. This allows
evaluation of all other hypergeometric functions of the same Schwarz type by 
using contiguous relations and Euler's formula (\ref{eq:euler}).
For example, (\ref{fptetra1}) and (\ref{fptetra1b})
represent solutions of the hypergeometric equation with local
exponent differences 1/2,1/3,1/3 at $z=1,0,\infty$, respectively;
and (\ref{fptetra1a}) and (\ref{fptetra1z}) evaluate contiguous
functions to (\ref{fptetra1}) and (\ref{fptetra1b}), respectively.

The formulas can be checked by expanding power series at $x=0$ on both sides,
also in the cases of genus 1 Darboux curves, as $\xi=O(\sqrt{x})$ for these curves.


\subsection{Tetrahedral hypergeometric equations}
\label{tetrahedral}

For the Schwarz type $(1/2,\,1/3,\,1/3)$ we give the following
evaluations. The Darboux covering is evident from the arguments of
the $\hpgo{2}{1}$ series.
\begin{eqnarray} \label{fptetra1}
\hpg{2}{1}{1/4,-1/12}{2/3}{\frac{x\,(x+4)^3}{4(2x-1)^3}} & =
&\big(1-2x\big)^{-1/4}.\\ 
\label{fptetra1a}\hpg{2}{1}{5/4,-1/12}{5/3}{\frac{x\,(x+4)^3}{4(2x-1)^3}}
& = & \frac{1+x}{\left(1+\frac14x\right)^2}\,\big(1-2x\big)^{-1/4}.\\
\label{fptetra1b}\hpg{2}{1}{1/4,\;7/12}{4/3}{\frac{x\,(x+4)^3}{4(2x-1)^3}}
& = & \frac{1}{1+\frac14x}\,\left(1-2x\right)^{3/4}.\\
\label{fptetra1z}\hpg{2}{1}{1/4,-5/12}{1/3}{\frac{x\,(x+4)^3}{4(2x-1)^3}}
& = & \big(1+{\textstyle\frac52}x\big)\,\big(1-2x\big)^{-5/4}.
\end{eqnarray}
For the Schwarz type $(1/3,\,1/3,\,2/3)$ we give the following
evaluations. 
\begin{eqnarray} \label{fptetra2}
\hpg{2}{1}{1/2,\,-1/6}{2/3}{\frac{x\,(x+2)^3}{(2x+1)^3}}
& = & \big(1+2x\big)^{-1/2}.\\
\label{fptetra2a}\hpg{2}{1}{1/2,\;5/6}{2/3}{\frac{x\,(x+2)^3}{(2x+1)^3}}
& = & \frac{1}{(1-x)^2}\;\big(1+2x\big)^{3/2}.\\ \label{fptetra2b}
\hpg{2}{1}{1/6,\;5/6}{4/3}{\frac{x\,(x+2)^3}{(2x+1)^3}} & = &
\frac{1}{1+\frac{1}{2}x}\,\big(1+2x\big)^{1/2}\,\big(1+x\big)^{1/3}. \\
\label{fptetra2z}\hpg{2}{1}{1/6,-1/6}{1/3}{\frac{x\,(x+2)^3}{(2x+1)^3}}
& = & \big(1+2x\big)^{-1/2}\,\big(1+x\big)^{1/3}.
\end{eqnarray}

\subsection{Octahedral hypergeometric equations}

For the Swartz type $(1/2,\,1/3,\,1/4)$ we have the following evaluations. 
\begin{eqnarray} \label{fpocta1}
\hpg{2}{1}{7/24,\,-1/24}{3/4}{\frac{108\,x\,(x-1)^4}{(x^2+14x+1)^3}}
&=& \big(1+14x+x^2\big)^{-1/8}.\\
\hpg{2}{1}{7/24,\,\;23/24}{7/4}{\frac{108\,x\,(x-1)^4}{(x^2+14x+1)^3}}
&=& \frac{1+2x-\frac{1}{11}x^2}{(1-x)^3}\,\big(1+14x+x^2\big)^{7/8}.\\
\hpg{2}{1}{5/24,\;13/24}{5/4}{\frac{108\,x\,(x-1)^4}{(x^2+14x+1)^3}}
&=& \frac{1}{1-x}\,\big(1+14x+x^2\big)^{5/8}.\\
\hpg{2}{1}{5/24,\,-11/24}{1/4}{\frac{108\,x\,(x-1)^4}{(x^2+14x+1)^3}}
&=& \frac{1-22x-11x^2}{\big(1+14x+x^2\big)^{11/8}}.
\end{eqnarray}
For the Swartz type
$(1/4,\,1/4,\,2/3)$ we have the following evaluations. 
\begin{eqnarray} \label{fpocta2}
\hpg{2}{1}{7/12,\,-1/12}{3/4}{\frac{27\,x\,(x+1)^4}{2(x^2+4x+1)^3}}
& = & \frac{\left(1+\frac12x\right)^{1/4}}{(1+4x+x^2)^{1/4}}.\\
\hpg{2}{1}{7/12,\;11/12}{7/4}{\frac{27\,x\,(x+1)^4}{2(x^2+4x+1)^3}}
& = & \frac{(1+\frac12x)^{1/4}\;(1+4x+x^2)^{7/4}}{(1+x)^3}.\\
\hpg{2}{1}{1/6,\;5/6}{5/4}{\frac{27\,x\,(x+1)^4}{2(x^2+4x+1)^3}}
& = & \frac{\big(1+2x\big)^{1/4}\,\big(1+4x+x^2\big)^{1/2}}{1+x}.\\
\hpg{2}{1}{1/6,\,-1/6}{1/4}{\frac{27\,x\,(x+1)^4}{2(x^2+4x+1)^3}}
& = & \frac{\big(1+2x\big)^{1/4}}{\big(1+4x+x^2\big)^{1/2}}.
\end{eqnarray}

\subsection{Icosahedral hypergeometric equations}

The Darboux covering for hypergeometric equations of the Schwarz
type $(1/2,1/3,1/5)$ is the following:
\begin{equation} \label{isophi1}
\varphi_1(x)=\frac{1728\;x\;(x^2-11x-1)^5}{(x^4+228x^3+494x^2-228x+1)^3}
\end{equation}
The simplest evaluations for this Schwarz type are these:
\begin{eqnarray} \label{fpicosa1} 
\hpg{2}{1}{19/60,\,-1/60}{4/5}{\varphi_1(x)} \equal 
\big(1-228x+494x^2+228x^3+x^4\big)^{-1/20}.\\ \label{fpicosa1a} 
\hpg{2}{1}{19/60,\,59/60}{4/5}{\varphi_1(x)} \equal
\frac{(1\!+\!66x\!-\!11x^2)(1\!-\!228x\!+\!494x^2\!+\!228x^3\!+\!x^4)^{19/20}}
{(1+x^2)\,(1+522x-10006x^2-522x^3+x^4)}.\\ 
\label{fpicosa1b}  
\hpg{2}{1}{11/60,\,31/60}{6/5}{\varphi_1(x)} \equal
\frac{(1\!-\!228x\!+\!494x^2\!+\!228x^3\!+\!x^4)^{11/20}}
{1+11x-x^2}.\\  \label{fpicosa1z}
\hpg{2}{1}{11/60,\,-29/60}{1/5}{\varphi_1(x)} \equal
\frac{1+435x-6670x^2-3335x^4-87x^5}
{(1\!-\!228x\!+\!494x^2\!+\!228x^3\!+\!x^4)^{29/20}}.
\end{eqnarray}
The Darboux covering for hypergeometric equations of the Schwarz
type $(1/2,1/3,2/5)$ is the same $\varphi_1(x)$. Here are the simplest contiguous
evaluations:
\begin{eqnarray}
\hpg{2}{1}{13/60,\,-7/60}{3/5}{\varphi_1(x)}  \equal 
\frac{1-7x}{\big(1-228x+494x^2+228x^3+x^4\big)^{7/20}}.\\
\hpg{2}{1}{13/60,\,53/60}{3/5}{\varphi_1(x)}  \equal 
\frac{(1\!+\!119x\!+\!187x^2\!+\!17x^3)(1\!-\!228x\!+\!494x^2\!+\!228x^3\!+\!x^4)^{13/20}}
{(1+x^2)\,(1+522x-10006x^2-522x^3+x^4)}.\nonumber\\
\hpg{2}{1}{17/60,\,37/60}{7/5}{\varphi_1(x)}  \equal 
\frac{\left(1+\frac{1}{7}x\right)(1\!-\!228x\!+\!494x^2\!+\!228x^3\!+\!x^4)^{17/20}}
{\left(1+11x-x^2\right)^2}.\\
\hpg{2}{1}{17/60,\,-23/60}{2/5}{\varphi_1(x)}  \equal 
\frac{(1\!+\!207x\!-\!391x^2\!+\!1173x^3\!+\!46x^4)}
{(1\!-\!228x\!+\!494x^2\!+\!228x^3\!+\!x^4)^{23/20}}.
\end{eqnarray}
The Darboux covering for the Schwarz type
$(1/2,1/5,2/5)$ is the following:
\begin{equation} \label{isophi2}
\varphi_2(x)= \frac{64\,x\,(x^2-x-1)^5}{(x^2-1)\,(x^2+4x-1)^5}.
\end{equation}
The simplest evaluations are the following:
\begin{eqnarray}
\hpg{2}{1}{7/20,\,-1/20}{4/5}{\varphi_2(x)} \equal 
\frac{(1+x)^{7/20}}{(1-x)^{1/20}\,(1-4x-x^2)^{1/4}}.\\
\hpg{2}{1}{7/20,\,19/20}{4/5}{\varphi_2(x)}  \equal 
\frac{(1+3x)\,(1+x)^{7/20}\,(1-x)^{19/20}\,(1-4x-x^2)^{7/4}}
{(1+x^2)\,(1+22x-6x^2-22x^3+x^4)}.\\
\hpg{2}{1}{3/20,\,11/20}{6/5}{\varphi_2(x)}  \equal 
\frac{(1+x)^{3/20}\,(1-x)^{11/20}\,(1-4x-x^2)^{3/4}}
{1+x-x^2}. \\
\hpg{2}{1}{3/20,\,-9/20}{1/5}{\varphi_2(x)}  \equal 
\frac{(1+12x-6x^2-2x^3-9x^4)\;(1+x)^{3/20}}{(1-x)^{9/20}\,(1-4x-x^2)^{9/4}}.
\end{eqnarray}
For other icosahedral Schwarz types, the Darboux curves are not
rational but have genus 1. We use then two algebraically
related variables $x,\xi$ in Algorithm \ref{kleinalg}.


The Darboux curve for hypergeometric equations of  the
Schwarz type $(1/3,1/3,2/5)$ is given by the equation
\begin{equation} \label{darbouxc1}
E_3: \qquad \xi^2=x\;(1+33x-9x^2).
\end{equation}
The Darboux covering is
\begin{equation} \label{isophi3}
\varphi_3(x,\xi)=\frac{144\;\xi\,\left(1+33x-9x^2\right)^2\,(1-9\xi+54x)}
{\left(1+21\xi-117x+9x\xi-234x^2\right)^3}
\end{equation}
Here are simplest evaluations for the Schwarz type $(1/3,1/3,2/5)$.
\begin{eqnarray} \label{icosellipta}
\hpg{2}{1}{3/10,\,-1/30}{3/5}{\varphi_3(x,\xi)}  \equal 
\frac{(1-9\xi+54x)^{1/30}}{(1+21\xi-117x+9x\xi-234x^2)^{1/10}}.\\
\hpg{2}{1}{3/10,\,29/30}{3/5}{\varphi_3(x,\xi)}  \equal  \nonumber \\
& & \hspace{-2.1cm} \label{icoselliptb3}
\frac{(1\!+\!21\xi\!-\!117x\!+\!9x\xi\!-\!234x^2)^{9/10}(1\!+\!9x)^2(1\!+\!198x\!-\!99x^2)}
{(1-9\xi+54x)^{29/30}\,(1\!-\!21\xi\!-\!117x\!-\!9x\xi\!-\!234x^2)^2}.\\
\hpg{2}{1}{7/10,\;11/30}{7/5}{\varphi_3(x,\xi)}  \equal 
\frac{\left(1+21\xi-117x+9x\xi-234x^2\right)^{11/10}}
{(1-9\xi+54x)^{11/30}\,(1+33x-9x^2)}. \label{icoselliptc3}\\
\hpg{2}{1}{\!-3/10,\,11/30}{2/5}{\varphi_3(x,\xi)}  \equal 
\frac{(1-9\xi+54x)^{19/30}\,(1-15\xi-72x-54x^2)}
{(1+21\xi-117x+9x\xi-234x^2)^{9/10}\,(1+9x)}. \label{icoselliptd3}
\end{eqnarray}
The Darboux curve and covering for the Schwarz type
$(1/3,2/3,1/5)$ are the same as in (\ref{darbouxc1}) and
(\ref{isophi3}). The simplest evaluations are:
\begin{eqnarray} \label{icoselliptk3}
\hpg{2}{1}{-1/10,\,17/30}{4/5}{\varphi_3(x,\xi)}  \equal  \frac{
(1-9\xi+54x)^{13/30}}{\left(1+21\xi-117x+9x\xi-234x^2\right)^{3/10}}.\\
\hpg{2}{1}{9/10,\,17/30}{9/5}{\varphi_3(x,\xi)}  \equal  \frac{
\left(1+21\xi-117x+9x\xi-234x^2\right)^{17/10}\left(1+\frac37x\right)}
{(1-9\xi+54x)^{17/30}\,\left(1+33x-9x^2\right)^2}.\label{icoselliptl3}\\
\hpg{2}{1}{1/10,\,23/30}{6/5}{\varphi_3(x,\xi)}  \equal  \frac{
\left(1\!+\!21\xi\!-\!117x\!+\!9x\xi\!-\!234x^2\right)^{3/10}\!
(1\!-\!9\xi\!+\!54x)^{7/30}(\xi\!+\!5x)}{\xi\;(1+9x)}.\label{icoselliptm3}\\
\label{icoselliptz3} \hpg{2}{1}{1/10,-7/30}{1/5}{\varphi_3(x,\xi)}
 \equal  \frac{ \!
(1-9\xi+54x)^{7/30}\,(1-21x)}{\left(1+21\xi-117x+9x\xi-234x^2\right)^{7/10}}.
\end{eqnarray}

The Darboux curve for the Schwarz type $(2/3,\,1/5,\,1/5)$ can be
given by
\begin{equation} \label{darbouxc2}
E_4: \qquad \xi^2=x\;(1+5x-5x^2).
\end{equation}
The Darboux covering is given by
\begin{equation} \label{isophi4}
\varphi_4(x,\xi)=\frac{432\,x\,\left(1-\frac{7}{5}\xi-9x-x^2\right)^5\,
(1\!+\!50x\!-\!125\xi^2\!+\!450x\xi\!-\!500x^2)}{(5\xi\!+\!57x)\,
\left(1\!+\!\frac{18}{5}\xi\!-\!16x\!+\!x^2\right)^5
(1\!+\!50x\!-\!125\xi^2\!-\!450x\xi\!-\!500x^2)}.
\end{equation}
The simplest evaluations are:
\begin{eqnarray} \label{icosellipta4}
\hpg{2}{1}{1/6,-1/30}{4/5}{\varphi_4(x,\xi)} \equal  \frac{
\left(1-\frac35\xi-\frac{34}{5}x\right)^{1/6}}
{(1\!+\!3\xi\!-\!20x)^{1/6}\left(
1\!+\!50x\!-\!125\xi^2\!-\!450x\xi\!-\!500x^2\right)^{1/30}}.\\
\hpg{2}{1}{1/6,29/30}{4/5}{\varphi_4(x,\xi)}  \equal  \nonumber \\
\label{icoselliptb4}& & \hspace{-3.4cm} \frac{(1+3\xi-20x)^{5/6}\,
\left(1-\frac35\xi-\frac{34}{5}x\right)^{1/6}\,
\left(1-\frac{35}4\xi-\frac{101}4x\right)}
{\left(1\!+\!50x\!-\!125\xi^2\!-\!450x\xi\!-\!500x^2\right)^{1/30}\!
\left(1\!-\!\frac{95}4\xi\!+\!\frac{83}4x\!+\!\frac{21}4\xi^2
\!+\!\frac{475}4x\xi\!+\!10x^2 \right)}.\\ \label{icoselliptc4}
\hpg{2}{1}{1/6,\,11/30}{6/5}{\varphi_4(x,\xi)}  \equal 
\left(1\!+\!50x\!-\!125\xi^2\!-\!450x\xi\!-\!500x^2\right)^{11/30}
\times \nonumber \\ && \hspace{-22pt}
\frac{(1+3\xi-20x)^{5/6}\left(1-\frac35\xi-\frac{34}5x\right)^{1/6}
\left(1+\frac{21}4\xi+\frac{41}4x\right)}
{(1-9x)\,\left(1-\frac74\xi-\frac{15}2x\right)\,(1+5\xi+10x)}.\\
\hpg{2}{1}{1/6,\,11/30}{1/5}{\varphi_4(x,\xi)}  \equal 
\left(1\!+\!50x\!-\!125\xi^2\!-\!450x\xi\!-\!500x^2\right)^{11/30}
\times \nonumber \\ && \hspace{-22pt}\label{icoselliptd4}
\frac{(1+3\xi-20x)^{5/6}\left(1-\frac35\xi-\frac{34}5x\right)^{1/6}
\left(1+\frac{21}4\xi+\frac{41}4x\right)}
{\left(1\!-\!\frac{95}4\xi\!+\!\frac{83}4x\!+\!\frac{21}4\xi^2
\!+\!\frac{475}4x\xi\!+\!10x^2 \right)\,(1+5\xi+10x)}.
\end{eqnarray}
The Darboux curve and covering for the Schwarz type
$(1/3,2/5,3/5)$ are the same as in (\ref{darbouxc2}) and
(\ref{isophi4}). The simplest evaluations are:
\begin{eqnarray}\label{icoselliptk4}
\hpg{2}{1}{\!-1/6,\,13/30}{3/5}{\varphi_4(x,\xi)} \equal 
\frac{ (1\!+\!50x\!-\!125\xi^2\!-\!450x\xi\!-\!500x^2)^{13/30}
\left(1-3\xi+2x\right)}{\left(1+3\xi-20x\right)^{5/6}\,
\left(1-\frac35\xi-\frac{34}{5}x\right)^{1/6}\,(1+5\xi+10x)}.\\
\hpg{2}{1}{\!5/6,\,13/30}{3/5}{\varphi_4(x,\xi)} \equal
(1\!+\!50x\!-\!125\xi^2\!-\!450x\xi\!-\!500x^2)^{13/30}\times\nonumber\\
\label{icoselliptl4}&&\hspace{-35pt}\frac{\left(1+3\xi-20x\right)^{13/6}
\left(1-\frac74\xi+\frac{25}2x-\frac{245}4x^2\right)
\left(1-\frac7{20}\xi-\frac{79}{20}x\right)}
{\left(1\!-\!\frac35\xi\!-\!\frac{34}{5}x\right)^{1/6}\!
\left(1\!-\!\frac{95}4\xi\!+\!\frac{83}4x\!+\!\frac{21}4\xi^2
\!+\!\frac{475}4x\xi\!+\!10x^2\right)^2\!(1-5x)^2}.\\
\hpg{2}{1}{5/6,\,7/30}{7/5}{\varphi_4(x,\xi)} \equal 
\left(1\!+\!50x\!-\!125\xi^2\!-\!450x\xi\!-\!500x^2\right)^{7/30}
\times \nonumber \\ && \hspace{-30pt} \label{icoselliptm4}
\frac{\left(1\!+\!\frac{18}{5}\xi\!-\!16x\!+\!x^2\right)^{7/6}
\!\left(1\!+\!\frac{1}{25}x\right)^{5/6}\!(1\!+\!5\xi\!+\!10x)}
{\left(1-\frac{7}{5}\xi-9x-x^2\right)^2\;(1-5x)^{7/6}}.\\
\label{icoselliptz4}
\hpg{2}{1}{\!-1/6,\,7/30}{2/5}{\varphi_4(x,\xi)} \equal
\frac{\left(1\!+\!50x\!-\!125\xi^2\!-\!450x\xi\!-\!500x^2\right)^{7/30}
\!\left(1\!-\!\frac{27}5\xi\!+\!\frac{58}5x\!-\!2x^2\right)}
{\left(1\!+\!\frac{18}{5}\xi\!-\!16x\!+\!x^2\right)^{5/6}
(1+\frac1{25}x)^{1/6}\,(1-5x)^{7/6}}.
\end{eqnarray}

The Darboux curve for the Schwarz type $(1/3,\,1/5,\,3/5)$ can be
given by
\begin{equation} \label{darbouxc3}
E_5: \qquad \xi^2=x\;(1+x)\,(1+16x).
\end{equation}
The Darboux covering can be given by
\begin{equation} \label{isophi5}
\varphi_5(x,\xi)=-\frac{54\;(\xi+5x)^3\,(1-2\xi+6x)^5}
{(1\!-\!16x^2)\,(\xi-5x)^2\,(1-2\xi-14x)^5}.
\end{equation}
\begin{eqnarray} \label{icosellipta5}
\hpg{2}{1}{\!-1/15,\,8/15}{4/5}{\varphi_5(x,\xi)} \equal
\frac{(1+4x)^{8/15}\,(\xi+5x)^{1/6}\;x^{1/15}}
{(1-2\xi-14x)^{1/3}\;(\xi-3x)^{3/10}}.\\
\hpg{2}{1}{14/15,\,8/15}{9/5}{\varphi_5(x,\xi)} \equal \nonumber \\
 & & \hspace{-60pt} \frac{(1\!-\!2\xi\!-\!14x)^{8/3}
\left(1\!+\!\frac23\xi\!+\!\frac23x\!-\!\frac{16}3x^2\right)
(1\!+\!4x)^{8/15}\,(\xi\!-\!5x)^2\,x^{1/15}}
{(1-2\xi+6x)^4\;(\xi+5x)^{11/6}\;(\xi-3x)^{3/10}}.\\
\hpg{2}{1}{2/15,\,11/15}{6/5}{\varphi_5(x,\xi)} \equal \frac{
(1\!-\!\xi\!+\!x)(1\!-\!2\xi\!-\!14x)^{2/3}(\xi\!+\!5x)^{1/6}(\xi\!-\!3x)^{13/10}}
{(1+\xi+x)\,(1-2\xi+6x)\,(1+4x)^{13/15}\;x^{11/15}}.\\
\label{icoselliptz5}
 \hpg{2}{1}{2/15,-4/15}{1/5}{\varphi_5(x,\xi)} \equal
\frac{(1\!+\!3\xi\!+\!x)\,(1\!+\!4x)^{2/15}\,(\xi\!+\!5x)^{7/6}\,(\xi\!-\!3x)^{3/10}}
{(1+\xi+x)\;(1-2\xi-14x)^{4/3}\;x^{11/15}}.
\end{eqnarray}

The Darboux curve for the Schwarz type $(1/5,\,1/5,\,4/5)$ can be
given by
\begin{equation} \label{darbouxc4}
E_6: \qquad \xi^2=x\;(1+x-x^2).
\end{equation}
The Darboux covering is given by
\begin{equation} \label{isophi6}
\varphi_6(x,\xi)=\frac{16\;\xi\,\left(1+x-x^2\right)^2\,(1-\xi)^2}
{(1+\xi+2x)\,(1+\xi-2x)^5}.
\end{equation}
The simplest evaluations are:
\begin{eqnarray} \label{icosellipta6}
\hpg{2}{1}{7/10,-1/10}{4/5}{\varphi_6(x,\xi)} \equal \frac{
(1-\xi+2x)^{1/15}\,(1-\xi)^{3/5}}{(1+\xi+2x)^{7/30}\,\sqrt{1+\xi-2x}}\\
\hpg{2}{1}{7/10,\,9/10}{9/5}{\varphi_6(x,\xi)} \equal \frac{
(1-\xi+2x)^{1/15}\,(1+\xi+2x)^{23/30}\,(1+\xi-2x)^{7/2}}{(1-\xi)^{7/5}\,(1+x-x^2)^2}\\
\hpg{2}{1}{1/10,\,9/10}{6/5}{\varphi_6(x,\xi)} \equal \frac{
(\xi+2x+x^2)\;(1+\xi)^{1/10}\,(1-\xi)^{3/10}}{\xi\,(1-\xi+2x)^{1/30}
(1+\xi+2x)^{2/15}\sqrt{1+\xi-2x}}.\\ \label{icosellipta6z}
\hpg{2}{1}{1/10,-1/10}{1/5}{\varphi_6(x,\xi)} \equal \frac{
(1+\xi)^{1/10}\,(1-\xi)^{3/10}}{(1-\xi+2x)^{1/30}\,(1+\xi+2x)^{2/15}\,\sqrt{1+\xi-2x}}.
\end{eqnarray}
The Darboux curve and covering for the Schwarz type
$(2/5,2/5,2/5)$ are the same as in (\ref{darbouxc4}) and
(\ref{isophi6}). The simplest evaluations are:
\begin{eqnarray}
\hpg{2}{1}{3/10,-1/10}{3/5}{\varphi_6(x,\xi)} \equal \frac{
(1-\xi+2x)^{2/15}\,(1+\xi+2x)^{1/30}\,(1-\xi)^{1/5}}{\sqrt{1+\xi-2x}}.\\
\hpg{2}{1}{3/10,\,9/10}{8/5}{\varphi_6(x,\xi)} \equal \nonumber \\
& & \hspace{-3cm} \frac{
(1\!-\!\xi\!+\!2x)^{2/15}(1\!+\!\xi\!+\!2x)^{1/30}(1\!-\!\xi)^{1/5}(1\!+\!\xi\!-\!2x)^{3/2}
\left(1\!+\!\frac{1}{2}\xi\!+\!\frac{1}{2}x\right)}{(1+x-x^2)\,(1-\xi+x-x^2)}.\\
\hpg{2}{1}{3/10,\,7/10}{7/5}{\varphi_6(x,\xi)} \equal \frac{
(1+\xi+2x)^{7/30}\,(1+\xi)^{1/5}\,(1+\xi-2x)^{3/2}}
{(1-\xi+2x)^{1/15}\,(1-\xi)^{2/5}\,(1+x-x^2)}.\\
\label{icoselliptz}\hpg{2}{1}{3/10,\,-3/10}{2/5}{\varphi_6(x,\xi)}
\equal \frac{(1+\xi+2x)^{7/30}\,(1+\xi)^{1/5}\,(1\!-\!3\xi\!+\!4x\!-\!2x^2)}
{(1-\xi+2x)^{1/15}\,(1-\xi)^{2/5}\,(1+\xi-2x)^{3/2}}.
\end{eqnarray}

\bibliographystyle{alpha}
\bibliography{../hypergeometric}

\begin{thebibliography}{BvHW03}

\bibitem[AAR99]{specfaar}
G.E. Andrews, R.~Askey, and R.~Roy.
\newblock {\em Special Functions}.
\newblock Cambridge Univ. Press, Cambridge, 1999.

\bibitem[BD79]{baltadw}
F.~Baldassarri and B.~Dwork.
\newblock On second order linear differential equations with algebraic
  solutions.
\newblock {\em American Journal of Mathematics}, 101:42--76, 1979.

\bibitem[Ber04]{maintphd}
M.~Berkenbosch.
\newblock {\em Algorithms and Moduli spaces for Differential Equations}.
\newblock PhD thesis, University of Groningen, 2004.

\bibitem[Beu02]{beukers}
F.~Beukers.
\newblock Gauss' hypergeometric function.
\newblock Technical report, Utrecht University, {\sf
  http://www.math.uu.nl/people/beukers/MRIcourse93.ps}, 2002.

\bibitem[Bou98]{boulang}
A.~Boulanger.
\newblock Contribution a l'etude des equations lineaires homogenes integrables
  algebriquement.
\newblock {\em Journal de l'Ecole Polytechnique}, 4:1--122, 1898.

\bibitem[Bri77]{brioschi}
F.~Brioschi.
\newblock La th\'eorie des formes dans l'int\'egration des \'equations
  diff\'erentielles line\'eaires du second ordre.
\newblock {\em Math. Annalen}, 11:401--411, 1877.

\bibitem[BvHW03]{WeHuBe}
M.~Berkenbosch, M.~van Hoeij, and J.-A. Weil.
\newblock Recent algorithms for solving second order differential equations.
\newblock Technical report, 2003.
\newblock Available at {\sf http://algo.inria.fr/seminars/sem01-02/weil.pdf}.

\bibitem[Fuc75]{fuchs}
L.~Fuchs.
\newblock Uber die linearen {D}ifferentialgleichungen zweiter {O}rdnung, welche
  algebraische {I}ntegralen besitzen, und eine neue {A}nwendung der
  {I}nvariantentheorie.
\newblock {\em Journ. fur die reine und angewandte {M}athematik}, 81:97--147,
  1875.

\bibitem[Kat72]{katz72}
N.~Katz.
\newblock Algebraic solutions of differential equations.
\newblock {\em Inv. Math.}, 18:1--118, 1972.

\bibitem[Kle77]{klein77}
F.~Klein.
\newblock Uber lineare differentialgleichungen {I}.
\newblock {\em Math. Annalen}, 11:115--118, 1877.

\bibitem[Kle78]{klein78}
F.~Klein.
\newblock Uber lineare differentialgleichungen {II}.
\newblock {\em Math. Annalen}, 12:167--179, 1878.

\bibitem[Kle84]{icosaklein}
F.~Klein.
\newblock {\em Vorlesungen \"uber das Ikosaeder und die Aufl\"osung del
  Gleichungen vom f\"unften Grade}.
\newblock Leipzig, 1884.

\bibitem[OY04]{ochiay}
H.~Ochiai and M.~Yoshida.
\newblock Polynomials associated with the contiguity relations of the
  hypergeometric functions with finite monodromy groups.
\newblock {\em International J. of Math.}, 15 (7):629--650, 2004.

\bibitem[Pep81]{pepin}
P.~Th. Pepin.
\newblock M\'ethodes pour obtenir les integrales algebriques des equations
  differentielles lineaires du second ordre.
\newblock {\em Atti dell'Accad. Pont. de Nouvi Lincei}, 36:243--388, 1881.

\bibitem[Sch72]{schwarz72}
H.A. Schwarz.
\newblock Ueber diejenigen {F}alle, in welchen die {G}aussische
  hypergeometrische {R}eihe eine algebraische {F}unktion ihres vierten
  {E}lements darstelt.
\newblock {\em Journ. f\"ur die reine und angewandte Math.}, 75:292--335, 1872.

\bibitem[SU93]{singulm2}
M.~F. Singer and F.~Ulmer.
\newblock Liouvillian and algebraic solutions of second and third order linear
  differential equations.
\newblock {\em Journ. Symb. Comp.}, 16(3):37--73, 1993.

\bibitem[vdPU98]{putulm}
M.~van~der Put and F.~Ulmer.
\newblock Differential equations and finite groups.
\newblock Technical report, Math. Sc. Res. Inst., Berkeley, California, 1998.

\bibitem[vHW05]{kleinvhw}
M.~van Hoeij and J.-A. Weil.
\newblock Solving second order linear differential equations with {K}lein's
  theorem.
\newblock Proceedings of the 2005 International Symposium on Symbolic and
  Algebraic Computation (ISSAC), 2005.

\bibitem[Vid03]{contiguous}
R.~Vid\=unas.
\newblock Contiguous relations of hypergeometric series.
\newblock {\em Journ.~Comp. Applied Math.}, 153:507--519, 2003.

\bibitem[Vid04]{algtgauss}
R.~Vid\=unas.
\newblock Algebraic transformations of {G}auss hypergeometric functions.
\newblock Available at {\sf http://arxiv.org/math.CA/0408269}, 2004.

\bibitem[Vid05a]{dalggaus}
R.~Vid\=unas.
\newblock Darboux evaluations of algebraic {G}auss hypergometric functions.
\newblock Available at {\sf http://arxiv.org/math.CA/0504264}, 2005.

\bibitem[Vid05b]{gammaeval}
R.~Vid\=unas.
\newblock Expressions for values of the gamma function.
\newblock {\em Kyushu J. Math.}, 59(2):267--283, 2005.

\bibitem[Vid08]{tdihedral}
R.~Vid\=unas.
\newblock Dihedral {G}auss hypergeometric functions.
\newblock Available at {\sf http://arxiv.org/abs/0807.4888}., 2008.

\end{thebibliography}

\end{document}